\documentstyle[12pt]{amsart}
\begin{document}
  \title{Triharmonic Riemannian submersions from 3-manifold of  constant curvature}
  \title[Triharmonic Riemannian submersion]
  {Triharmonic Riemannian submersions from 3-dimensional manifolds of  constant curvature}
 \author{Tomoya Miura}
  \address{Department of Mathematics, Shimane University, Matsue, 690-8504, Japan.}
  \email{tomo.miura3@@gmail.com}
  \author{Shun Maeta}
  \address{Department of Mathematics, Shimane University, Matsue, 690-8504, Japan.}
  \email{shun.maeta@@gmail.com and maeta@@riko.shimane-u.ac.jp}

    \keywords{harmonic maps, biharmonic maps, triharmonic maps, Riemannian submersions}
  \subjclass[2010]{primary 58E20, secondary 53C43.}
\thanks{The second author is partially supported by the Grant-in-Aid for Young Scientists(B), No.15K17542, Japan Society for the Promotion of Science.
}

\maketitle
\begin{abstract}
For biharmonic maps, there is a famous conjecture named Chen's conjecture.
In later paper, Wang and Ou gave an affirmative partial answer to submersion version of Chen's conjecture.
In this paper, we give an affirmative partial answer to submersion version of generalized Chen's conjecture, that is, triharmonic Riemannian submersions from a 3-dimensional space form into a surface is harmonic. 
    \end{abstract}
\numberwithin{equation}{section}
\theoremstyle{plain}
\newtheorem{df}{Definition}[section]
\newtheorem{th}[df]{Theorem}
\newtheorem{prop}[df]{Proposition}
\newtheorem{lem}[df]{Lemma}
\newtheorem{cor}[df]{Corollary}
\newtheorem{rem}[df]{Remark}
\newtheorem{cnj}[df]{Conjecture}

\section{Introduction}

A harmonic map which is a critical point of the energy functional\\ $E(\varphi)=\frac{1}{2}\int_M|d\varphi|^2v_g$ is important notion in geometry.

In 1964, Eells and Sampson \cite{ES} extended a harmonic map to a polyharmonic map of order $k$.
A polyharmonic map of order $k$ is defined a critical point of the $k$-energy functional (cf. \cite{EL2}, \cite{ES}) \\
\begin{equation}
E_k(\varphi)=\frac{1}{2}\int_M|(d+\delta)^k\varphi|^2v_g.
\label{eq:111}
\end{equation}
For $k=2$, the map called a biharmonic map and for $k=3$, the map called a triharmonic map.

A fundamental study of biharmonic maps is to classify all biharmonic maps between Riemannian manifolds.
For a biharmonic map, there is a famous conjecture.

\vspace{10pt}

{\bf Chen's conjecture} (cf. {\cite{CC}}): {\em Any biharmonic isometric immersion $\varphi: (M^m,g) \to {\mathbb{R}}^n$ into the Euclidean space is harmonic.} 

\vspace{10pt}

There are several affirmative partial answers to this conjecture.
In particular, for $M^2$ in ${\mathbb{E}}^3$ \cite{CC} and for $M^3$ in ${\mathbb{E}}^4$ \cite{DF}, \cite{HV}, this conjecture was proved.
Also, for $\varphi$ is a curve \cite{DI}  and $\varphi$ is proper \cite{SA1}, this conjecture was proved (see also \cite{sm12}, \cite{Luo2014}).

For the Chen's conjecture, the second author posed the generalized Chen's conjecture.

\vspace{10pt}

{\bf Maeta's Generalized Chen's conjecture} (cf. {\cite{SM1}, see also \cite{CC1}}):
 {\em The only $k$-harmonic submanifolds (polyharmonic submanifolds of order $k$) in Euclidean spaces are the
minimal ones.}

\vspace{10pt}

In later paper, Wang and Ou \cite{WO} showed the following.
\begin{th}[\cite{WO}]\label{Wang-Ou}
Let $\varphi : (M^3(c),g) \to (N^2,h)$ be a biharmonic Riemannian submersion from a space form of constant sectional curvature c. Then, $\varphi$ is harmonic.
\end{th}

We pose submersion version of generalized Chen's conjecture. 

\begin{cnj}
~\\
$1$. Are biharmonic submersions  $\varphi:M^m(c) \to N^{m-1}$ harmonic?\\
$2$. Are triharmonic submersions $\varphi:M^m(c) \to N^{m-1}$ harmonic?
\end{cnj}

In this paper, we give an affirmative partial answer to submersion version of generalized Chen's conjecture.
\begin{th}\label{main Th}
Let $\varphi:(M^3(c),g) \to (N^2,h)$ be a triharmonic Riemannian submersion from a space form of constant sectional curvature c. Then, $\varphi$ is harmonic.
\end{th}

\begin{rem}\label{bitri}
When we extend biharmonic immersions to triharmonic immersions, we have to add some additional assumptions (see~\cite{SM}).
However, in the above theorem, we do not need any additional assumptions.
\end{rem}
From the definition of a harmonic morphism (cf.~\cite{BW}), we have the following.
  
\begin{cor}
Let $\varphi:(M^3(c),g) \to (N^2,h)$ be a triharmonic Riemannian submersion from a space form of constant sectional curvature c. Then, $\varphi$ is harmonic morphism.

\end{cor}
By the same argument as in \cite{WO}, we obtain the following result.
\begin{cor}
$(1)$ If $\varphi: {\mathbb{R}}^3 \to (N^2,h)$ is a triharmonic horizontally homothetic submersion from Euclidean space, then $(N^2,h)$ is flat and $\varphi$ is a composition of an orthogonal projection ${\mathbb{R}}^3 \to {\mathbb{R}}^2$ followed by a covering map ${\mathbb{R}}^2 \to (N^2,h)$;\\
$(2)$ There exists no triharmonic Riemannian submersion $\varphi: {\mathbb{H}}^3 \to (N^2,h)$ no matter what $(N^2,h)$ is. 
\end{cor}
\section{Preliminaries}
In this section, to proof the main theorem, we prepare several definitions and notions. In this paper, all manifolds, maps and tensor fields are assumed to be smooth.
$\nabla$ and ${\nabla}^N$ denote the Levi-Civita connections on $(M,g)$ and $(N,h)$, respectively. $\nabla^{\varphi}$ denotes the induced connection on $\varphi^{-1}TN$.
The first variation formula of the trienergy is given as follows:
\begin{center}
$\left. \frac{d}{dt} \right|_{t=0} E_3({\varphi}_t) = -\int_M \langle {\tau}_3(\varphi), V \rangle v_g,$ \quad ($V\in \Gamma(\varphi^{-1}TN)$),
\end{center}
where,
\begin{align*}
\tau_3(\varphi)&=J(\overline{\triangle}(\tau(\varphi))) - \sum_{i=1}
^m R^N({\nabla}^{\varphi}_{e_i}\tau(\varphi),\tau(\varphi))d{\varphi}(e_i),\\
 J(\tau(\varphi))&=\overline{\triangle}(\tau(\varphi)) - \sum_{i=1}
^m R^N(\tau(\varphi),d{\varphi}(e_i))d{\varphi}(e_i),\\
\overline{\triangle}\tau(\varphi)&= \sum_{i=1}
^m {\nabla}^{\varphi}_{e_i}{\nabla}^{\varphi}_{e_i}\tau(\varphi) - {\nabla}^{\varphi}_{\nabla_{e_i}{e_i}}\tau(\varphi),\\
\intertext{and}
R^N(U,W)&= \nabla^N_U\nabla^N_W -\nabla^N_W\nabla^N_U-\nabla^N_{[U,W]}. \quad(U,W \in \Gamma(TN)). 
\end{align*}

$\tau^3(\varphi)$ is called the {\em tritension field} of $\varphi$.
A smooth map $\varphi$ of $(M^m,g)$ to $(N^n,h)$ is said to be $triharmonic$ if $\tau^3(\varphi)=0.$

We introduce main tools of Riemannian submersions from $3$-manifolds named adapted to a Riemannian submersion and integrability data (cf. {\cite{WO}}).
Let $\varphi : (M^3,g) \to (N^2,h)$ be a Riemannian submersion.
We can choose the special orthonormal frame adapted to $\varphi$ $\{e_1, e_2, e_3\}$ on $(M^3,g)$, where $e_3$ is vertical.  The frame called adapted to Riemannian submersion.
As is well known, the frame always exists. Since $e_3$ is vertical, $[e_1,e_3]$ and $[e_2,e_3]$ are vertical. Let $\{\epsilon_1,\epsilon_2\}$ be the orthonormal frame on $(N^2,h)$.
If we assume that 
\begin{center}
$[\epsilon_1,\epsilon_2] = F_1\epsilon_1 + F_2\epsilon_2$, \quad $(F_1, F_2 \in C^{\infty}(N))$,
\end{center}
then, we can assume
\begin{eqnarray}
  \left\{
    \begin{array}{l}
    [e_1, e_3] = k_1e_3,\\

    [e_2, e_3] = k_2e_3,\\      

    [e_1,e_2] = f_1e_1 +f_2e_2 -2{\sigma}e_3, 
   \end{array}
  \right.
\label{eq:11}
\end{eqnarray}
where,  $k_1,k_2,\sigma, f_1 = F_1\circ\varphi$ and $f_2 = F_2\circ\varphi \in C^{\infty}(M)$ called integrability data of the adapted to Riemannian submersion $\varphi$.

\section{Proof of the main theorem}

\begin{prop}
Let $\varphi :(M^3,g) \to (N^2,h)$ be a Riemannian submersion with the adapted frame $\{e_1,e_2,e_3\}$ and the integrability data $f_1,f_2,k_1,k_2$ and $\sigma$.
Then, the Riemannian submersion $\varphi$ is $triharmonic$ if and only if
\begin{eqnarray}
  \left\{
    \begin{array}{l}
      {\triangle}^MA+f_1e_1(B)+e_1(Bf_1)+f_2e_2(B)+e_2(Bf_2)-f_1k_1B-k_2f_2B\medskip \\
-A\{K^N+f^2_1+f^2_2\}+K^N(-e_2(k_1)k_2+e_2(k_2)k_1-k^2_1f_2-k^2_2f_2)=0,\bigskip \\

      {\triangle}^MB-f_1e_1A-e_1(Af_1)-f_2e_2(A)-e_2(Af_2)+f_1k_1A+f_2k_2A \medskip\\
-B\{K^N+f^2_1+f^2_2\}+K^N(e_1(k_1)k_2-e_1(k_2)k_1+f_1k^2_2+f_1k^2_1)=0,
    \end{array}
  \right.
\label{eq:a}
\end{eqnarray}

where,
\begin{align*}
A&=-{\triangle}^Mk_1-f_1e_1(k_2)-e_1(k_2f_1)-f_2e_2(k_2)-e_2(k_2f_2)+k_1k_2f_1+f_2k^2_2\\
&+k_1f^2_1+k_1f^2_2,
\end{align*}
\begin{align*}
B&=-{\triangle}^Mk_2+f_1e_1(k_1)+e_1(k_1f_1)+f_2e_2(k_1)+e_2(k_1f_2)-k_1k_2f_2-f_1k^2_1\\
&+k_2f^2_1+k_2f^2_2,
\end{align*}
${\triangle}^M(k_1)=\displaystyle \sum_{i=1}^3 e_ie_i(k_1)+f_1e_2(k_1)-f_2e_1(k_1)-k_1e_1(k_1)-k_2e_2(k_1),$
\\and\\
$K^N = R^N_{1212}\circ\pi= e_1(f_2)-e_2(f_1)-f^2_1-f^2_2$.

\end{prop}
{\it Proof}
A simple computation, using (\ref{eq:11}) and Koszul formula, yields
\begin{align*}
& {\nabla}_{e_1}e_1 = -f_1e_2,                  {\nabla}_{e_1}e_2 = f_1e_1-{\sigma}e_3,    {\nabla}_{e_1}e_3 = {\sigma}e_2,\\ 
 &{\nabla}_{e_2}e_1 = -f_2e_2+ {\sigma}e_3,   {\nabla}_{e_2}e_2 = f_2e_1,                      {\nabla}_{e_2}e_3 = -{\sigma}e_1,\\
 &{\nabla}_{e_3}e_1 = -k_1e_3+ {\sigma}e_2,   {\nabla}_{e_3}e_2 = -{\sigma}e_1-k_2e_3,   {\nabla}_{e_3}e_3 = k_1e_1+k_2e_2.
\end{align*}
Furthermore, noting that $e_3$ is vertical, the following holds.
\begin{align*}
& {\nabla}^{\varphi}_{e_1}\epsilon_1 = -f_1{\epsilon}_2, \quad {\nabla}^{\varphi}_{e_1}\epsilon_2 = f_1{\epsilon}_1,   \\
& {\nabla}^{\varphi}_{e_2}\epsilon_1 = -f_2{\epsilon}_2, \quad {\nabla}^{\varphi}_{e_2}\epsilon_2 = f_2{\epsilon}_1,  \\    
 &{\nabla}^{\varphi}_{e_3}\epsilon_j = 0 \quad (j=1,2). 
\end{align*} 
Consequently, the tension field of ${\varphi}$ is given by
\begin{center}
$\tau(\varphi)=-k_1{\epsilon}_1-k_2{\epsilon}_2$.
\end{center}
Under these preparations, we calculate ${\tau}^3(\varphi)=0.$
First, by a direct computation, we get
\begin{center}
$\overline{\triangle}\tau(\varphi) = A{\epsilon}_1 + B{\epsilon}_2$.
\end{center}
Second, a straightforward computation  yields

\begin{align*}
J(\overline{\triangle}\tau(\varphi)) &= [{\triangle}^MA+f_1e_1(B)+e_1(Bf_1)+f_2e_2(B)+e_2(Bf_2)\\ 
&-f_1k_1B-k_2f_2B-A(f^2_1+f^2_2)-AK^N]{\epsilon}_1\\
&+[{\triangle}^MB-f_1e_1(A)-e_1(Af_1)-f_2e_2(A)-e_2(Af_2)\\
&+f_1k_1A+f_2k_2A-B(f^2_1+f^2_2)-BK^N]{\epsilon}_2.
\end{align*}
Third, in the same way, we get
\begin{align*}
-R^N({\nabla}^{\varphi}_{e_i}\tau(\varphi),\tau(\varphi))d{\varphi}(e_i)&=K^N(-e_2(k_1)k_2+e_2(k_2)k_1-k^2_1f_2-k^2_2f_2){\epsilon}_1\\
&+K^N(e_1(k_1)k_2-e_1(k_2)k_1+k^2_1f_1+k^2_2f_1){\epsilon}_2.
\end{align*}
Therefore, we can get (\ref{eq:a}).
\qed

\vspace{10pt}

If $M^3$ has a constant sectional curvature $c$, for all orthonormal frame $\{e_1,e_2,e_3\}$ on $M^3$ adapted to Riemannian submersion and for all integrability date  $f_1,f_2,k_1,k_2$ and $\sigma$, we have two facts. (cf.~\cite{WO}):
\begin{equation}
e_3(f_1)=e_3(f_2)=e_3(k_1)=e_3(k_2)=e_3(\sigma)=0,
\label{eq:kk}
\end{equation}
and we can choose adapted orthonormal frame $\{e_1,e_2,e_3\}$ such that $k_2=0.$

These facts simplify the triharmonic equation.

\begin{cor}\label{cor2}
Let $\varphi :(M^3(c),g) \to (N^2,h)$ be a Riemannian submersion from a space form of constant sectional curvature c with an adapted frame $\{e_1,e_2,e_3\}$ and the integrability date $f_1,f_2,k_1,k_2$ and $\sigma$ with $k_2=0$. Then, the Riemannian submersion $\varphi$ is triharmonic if and only if
\begin{eqnarray}
  \left\{
    \begin{array}{l}
      {\triangle}^MA'+f_1e_1(B')+e_1(B'f_1)+f_2e_2(B')+e_2(B'f_2)-f_1k_1B' \\
-A'\{K^N+f^2_1+f^2_2\}-K^N(k^2_1f_2)=0 ,\bigskip\\
      {\triangle}^MB'-f_1e_1A'-e_1(A'f_1)-f_2e_2(A')-e_2(A'f_2)+f_1k_1A' \\
-B'\{K^N+f^2_1+f^2_2\}+K^N(k^2_1f_1)=0,
    \end{array}
  \right.
\end{eqnarray}
where,
\\
$A'=-{\triangle}^Mk_1+k_1f^2_1+k_1f^2_2$,
\\and\\
$B'=f_1e_1(k_1)+e_1(k_1f_1)+f_2e_2(k_1)+e_2(k_1f_2)-f_1k^2_1$.

\end{cor}

By using Corollary $\ref{cor2},$ we prove the main theorem.
\begin{th}(Theorem $1.2.$)
Let $\varphi:(M^3(c),g) \to (N^2,h)$ be a triharmonic Riemannian submersion from a space form of constant sectional curvature c. Then, $\varphi$ is harmonic.
\end{th}

{\it Proof} \quad 
By the facts, we choose an orthonormal frame $\{e_1,e_2,e_3\}$ adapted to the Riemannian submersion
 with integrability data $f_1,f_2,k_1,k_2$ and $\sigma$
such that $k_2=0.$
Since $M^3(c)$ has constant sectional curvature c,
the following holds.
\begin{eqnarray}
  \left\{
    \begin{array}{l}
      R^M_{1312}=e_1(\sigma)-2k_1\sigma=0.\\
      R^M_{1313}=-[-e_1(k_1)-\sigma^2+k^2_1]=c.\\
      R^M_{1323}=k_1f_1=0.\\
      R^M_{1212}=-(e_2(f_1)-e_1(f_2)+f^2_1+f^2_2+3\sigma^2)=c.\\
     R^M_{1223}=e_2(\sigma)=0.\\
     R^M_{2313}=e_2(k_1)=0.\\
     R^M_{2323}=-(-\sigma^2+k_1f_2)=c.\\ 
    \end{array}
  \right.
\label{eq:b}
\end{eqnarray}
Here, by the 3rd equation in (\ref{eq:b}), if $k_1=0$ on $M^3(c)$, $\varphi$ is harmonic. 
Therefore, we consider $k_1\neq0$ on some set $\Omega \subset M^3(c)$ (in this case, $f_1=0$ on $\Omega$) and prove that the case cannot happen.

By the 1st equation of  (\ref{eq:b}), we get
\begin{equation}
 e_1(\sigma)=2k_1\sigma.
\label{eq:s}
\end{equation}
By the 2nd and the 7th equations in (\ref{eq:b}), we get
\begin{equation}
e_1(k_1)=k^2_1-k_1f_2.
\label{eq:k}
\end{equation}
By applying $e_1$ to both sides of the 7th equation in (\ref{eq:b}), we get
\begin{equation}
e_1(f_2)= f^2_2-k_1f_2+4\sigma^2.
\label{eq:ff}
\end{equation}
Also, in this case, we get 
\begin{equation}
A'=-k^3_1+2f_2k^2_1+4\sigma^2k_1.
\label{eq:aa}
\end{equation}
By using (\ref{eq:s}), (\ref{eq:k}), (\ref{eq:ff}) and (\ref{eq:aa}), the triharmonic equation is given as

\begin{equation}
\triangle^MA'-A'(e_1(f_2))-(e_1(f_2)-f^2_2)(f_2k^2_1)=0.
\label{eq:0}
\end{equation}
Here, by (\ref{eq:b}) and (\ref{eq:kk}), we get
\begin{center}
 $e_2(A')=e_3(A')=0$.
\end{center}
Consequently, laplacian of $A'$ can be rewritten as
\begin{center}
 ${\triangle}^M(A')=e_1e_1(A')-f_2e_1(A')-k_1e_1(A')$.
\end{center}
Since
\begin{align*}
e_1(A')&=-3k^4_1+5f_2k^3_1+(-2f^2_2+28\sigma^2)k^2_1-4\sigma^2f_2k_1,\\
\intertext{and}
e_1e_1(A')&=-12k^5_1+22f_2k^4_1+(-10f^2_2+188\sigma^2)k^3_1-88\sigma^2f_2k^2_1-16\sigma^4k_1,
\end{align*}
the equation of (\ref{eq:0}) can be written as
\begin{equation}
-9k^5_1+19f_2k^4_1+(164\sigma^2-9f^2_2)k^3_1-120\sigma^2f_2k^2_1-32\sigma^4k=0.
\label{eq:99}
\end{equation}
Since $k_1\neq0$ on $\Omega$, (\ref{eq:99}) implies
\begin{equation}
-9k^4_1+19f_2k^3_1+(164\sigma^2-9f^2_2)k^2_1-120\sigma^2f_2k_1-32\sigma^4=0.
\label{eq:c}
\end{equation}
Here, applying $e_1$ to both sides of (\ref{eq:c}), we get

\begin{equation}
-18k^4+37f_2k^3_1+(-19f^2_2+530\sigma^2)k^2_1-440\sigma^2f_2k_1-368\sigma^4=0.
\label{eq:d}
\end{equation}
By  $(\ref{eq:c})$ and $(\ref{eq:d})$, we get
\begin{equation}
f_2k^3+(-202\sigma^2+f^2_2)k^2_1+200{\sigma}^2f_2k_1+304\sigma^4=0.
\label{eq:e}
\end{equation}
Applying $e_1$ to both sides of (\ref{eq:e}) again, we get
\begin{equation}
f_2k^3_1+(-f^2_2-604\sigma^2)k^2_1+606\sigma^2f_2k_1+1616\sigma^4=0.
\label{eq:f}
\end{equation}
By (\ref{eq:e}) and (\ref{eq:f}), we get
\begin{equation}
(201\sigma^2+f^2_2)k^2_1-203\sigma^2f_2k_1-656\sigma^4=0.
\label{eq:g}
\end{equation}
By the same argument, we get
\begin{equation}
-f^2_2k^2_1+2\sigma^2f_2k_1-354\sigma^4=0.
\label{eq:ccc}
\end{equation}
Substituting the 7th equation of (\ref{eq:b}) into (\ref{eq:ccc}), we get
\begin{equation}
-353\sigma^4=c^2.
\label{eq:h}
\end{equation}
Applying $e_1$ to both sides of (\ref{eq:h}), we get
$$\sigma^4k_1=0.$$
From this, we get
$$\sigma=0~\text{on}~\Omega.$$
Substituting $\sigma=0$ into (\ref{eq:g}), we get 
\begin{center}
$k_1f_2=0$.
\end{center}
From the above arguments, we obtain $f_2=0$ on $\Omega$. Then, by using the 7th equation in  (\ref{eq:b}), we get $c=0$ on $\Omega$.
Thus, we get $e_1(k_1)=k^2_1.$
Hence, the triharmonic equation is given as\\
\begin{equation}
{\triangle}^Mk^3_1=0.
\label{eq:12}
\end{equation} 
Thus, it is easy to see that (\ref{eq:12}) is given by
$$k^5_1=0~\text{on}~\Omega.$$
Therefore, we get $k_1=0$ on $\Omega$ which is a contradiction.
Therefore, we complete the proof of the theorem.
\qed

\section{Appendix}

By the same argument as in the last section, we will give another proof of Theorem $\ref{Wang-Ou}$.
We obtain necessary and sufficient condition for a Riemannian submersion $\varphi:(M^3(c),g)\rightarrow (N^2,h)$ to be biharmonic as follows (cf.~\cite{WO}):

\begin{cor}\label{cor2-a}
Let $\varphi :(M^3(c),g) \to (N^2,h)$ be a Riemannian submersion from a space form of constant sectional curvature c with an adapted frame $\{e_1,e_2,e_3\}$ and the integrability date $f_1,f_2,k_1,k_2$ and $\sigma$ with $k_2=0$. Then, the Riemannian submersion $\varphi$ is biharmonic if and only if
\begin{eqnarray}
  \left\{
    \begin{array}{l}
      -\Delta ^M k_1+k_1(-K^N+f_1^2+f_2^2)= 0 ,\bigskip\\
      f_1e_1(k_1)+e_1(k_1f_1)+f_2e_2(k_1)+e_2(k_1f_2)-k_1^2f_1=0.
    \end{array}
  \right.
\end{eqnarray}
\end{cor}

\quad\\

{\it Proof of Theorem \ref{Wang-Ou}} \quad 
By the same argument as in the proof of Theorem $\ref{main Th}$, we obtain $(\ref{eq:b})$.
By the 3rd equation in (\ref{eq:b}),  if $k_1=0$ on $M^3(c)$, $\varphi$ is harmonic. 
Therefore, we consider $k_1\neq0$ on some set $\Omega \subset M^3(c)$ (in this case, $f_1=0$ on $\Omega$) and prove that the case cannot happen.

By the 1st equation of  (\ref{eq:b}), we get
\begin{equation}
 e_1(\sigma)=2k_1\sigma.
\label{eq:s-a}
\end{equation}
By the 2nd and the 7th equations in (\ref{eq:b}), we get
\begin{equation*}
e_1(k_1)=k^2_1-k_1f_2.
\label{eq:k-a}
\end{equation*}
By applying $e_1$ to both sides of the 7th equation in (\ref{eq:b}), we get
\begin{equation*}
e_1(f_2)= f^2_2-k_1f_2+4\sigma^2.
\label{eq:ff-a}
\end{equation*}
Since 
\begin{equation}\label{dk-a}
-\Delta^Mk_1=5k_1\sigma^2-k_1^3-k_1c+f_2(k_1^2-\sigma^2+c),
\end{equation}
the biharmonic equation is given as
\begin{equation}
k_1^2(3\sigma^2-k_1^2-3c)=0.
\label{eq:0-a}
\end{equation}
By the assumption $k_1\not=0$ on $\Omega$, we can get
\begin{equation}\label{a-1-a}
3\sigma^2-k_1^2-3c=0.
\end{equation}
Applying $e_1$ to this equation, we have
\begin{equation}\label{a-2-a}
k_1^2=7\sigma^2-c,
\end{equation}
since $k_1\not=0$ on $\Omega$.
Combining $(\ref{a-1-a})$, $(\ref{a-2-a})$ and $(\ref{eq:s-a})$, we have $k_1=\sigma=c=0$ on $\Omega$, which is a contradiction.
Therefore, we complete the proof of the theorem.
\qed


\vskip0.6cm\par


\end{document}